\documentclass[a4paper,11pt]{article}
\usepackage[english]{babel}
\usepackage{amsmath,amsthm,amsfonts,amssymb}
\usepackage{graphicx}
\usepackage{cite}
  \begin{document}    


\title{Numerical solution of the Volterra equations \\
of the first kind that appear in an inverse\\
 boundary-value problem of heat conduction
\thanks{This work was supported by RFBR grant № 12-01-00722-a}}

\author{Svetlana V. Solodusha\\
(Melentiev Energy Systems Institute SB RAS, Russia),\\
Natalia M. Yaparova\\
(South Ural State University, Chelyabinsk, Russia)}

\date{}
\maketitle

\abstract{The paper considers the integral Volterra equations of the first kind which are related to the inverse boundary-value heat conduction problem. The algorithms have been developed to numerically solve the respective integral equations, which are based on the midpoint rule and product integration method.  
 }
 
\section{	Problem statement}

Let us consider the following inverse boundary-value problem 
\begin{gather} \label{solod1}
u_{t} = u_{xx} ,\quad x \in \left( {0,1} \right), \quad t \geqslant
0,
\end{gather}
with boundary conditions 
\begin{gather} \label{solod2}
 u\left( {x,0} \right) = 0, \; u\left( {0,t} \right) = 0, \; u_{x} \left( {0,t} \right) = g\left( {t} \right),\;
x \in \left( {0,1} \right), \quad t \geqslant 0.
\end{gather}
It is necessary to find a boundary value of the function 
\begin{gather} \label{solod3}
u\left( {1,t} \right)=\phi(t), \quad t \geqslant 0.
\end{gather}

Let $g \in C^{2+\eta}_{[0,T]}$ at $T>0$, $\eta\in (0,1)$ and there exist constants   $M>0,$ $m\geq 0$ such that   $|g(t)|\leq M e^{mt}$ for $t\in
[0,T]$. 

 At  $g\left( {t} \right) = g_{0} \left( {t} \right)
$ there exists an accurate solution to  $u(1,t)=\phi_0 \left( {t} \right) $. Let us consider the case where instead of  $g_{0}(t) $ we know some approximations of  $g_{\delta }(t) $ and an error level   $\delta > 0$ such that  ${\left\| {g_{\delta}  - g_{0}} \right\|}_C \leqslant
\delta $. The uniqueness of the solution to problem (\ref{solod1}) --- (\ref{solod3}) is proved in  \cite{lavrentiev1980}.
 
To solve  (\ref{solod1}) --- (\ref{solod3}) we use the approach from  \cite{yaparova2013}, based on the direct and inverse Laplace transform, supposing that there exist constants  $C >0$ and $\beta _{0} \geqslant 0$ for which $\left| {\phi
\left( {t} \right)} \right| \leqslant C e^{\beta _{0} t}$ holds for $\forall x\in[0,1]$, 
$t \in \left[ {0,\infty}  \right)$ and  $\phi \left( {t} \right)$ meets the Dirichlet conditions $\forall\, t \in \left[ {0,T} \right]$. The author of  \cite{yaparova2013} suggests an algorithm to reduce (\ref{solod1}) --- (\ref{solod3})  to the integral Volterra equation of the first kind  
\begin{gather}\label{solod4}
 {A}\phi = 2\int\limits_{0}^{t} {\pi ^{2} \sum\limits_{p = 1}^{N} {\left( {
- 1} \right)^{p + 1} p^{2}e^{ - \pi ^2 p^2\left( {t - s}
\right)}}\phi \left( {s} \right)ds}  =g_{\delta}(t),\,\,
\end{gather}
\noindent
where $ 0 \leq s \leq t\leq T.$

\section{An algorithm for solving equation (\ref{solod4})}

Rewrite (\ref{solod4}) in the form  
\begin{gather}\label{solod5}
 {A}\phi = \int\limits_{0}^{t}K_N(t-s)\phi ( s)ds  =y(t),\,\, 0 \leq s \leq t\leq T,
\end{gather}
where
\begin{gather}\label{solod6}
 K_N(t-s)=\pi ^{2} \sum\limits_{p = 1}^{N} {\left( {
- 1} \right)^{p + 1} p^{2}e^{ - \pi ^2 p^2\left( {t - s}
\right)}}, \; y(t)=\frac{1}{2}g_{\delta}(t).
\end{gather}
To understand the specific features of the integral Volterra equation of the first kind  (\ref{solod5}), (\ref{solod6}) it is useful to consider the Volterra kernels  $K_{N} \in C_{\Delta}$, $\Delta=\{t,s/0 \leq s \leq t\leq T\}$ at fixed values of  $N$. Table 1 presents the values of  $K_N$ for $t=0$,  as well as the roots  $t^*$, that were obtained by solving the equations $K_N(t)=0$, $N=\overline{10,21}$ \cite{solodusha2014}. Figure 1 demonstrates the calculated values of  $t^{*}$ for  $N=\overline{10,41}$. 
\begin{table}[htb]
\caption{  Numerical characteristics of the Volterra kernels  $K_{N}$.}
\vspace{0.2cm}
\begin{center}
\begin{tabular}{|c|c|c|c|c|c|}\hline
$N$&$t^*$&$K_N(0)$&$N$&$t^*$&$K_N(0)$ \\  \hline
 10   &  0.01378 &-542.828& 16&0.00913&-1342.266 \\
 11   &  0.01221&651.394&17&0.00631&1510.049 \\
 12   & 0.01173&-769.829&18&0.00809&-1687.702 \\
 13 &  0.01022&898.134&19&0.00516&1875.225 \\
 14   &  0.01019&-1036.308&20&0.00735&-2072.617 \\
 15   &  0.00789&1184.353&21&0.00429&2279.879 \\   \hline
\end{tabular}
\end{center}
\end{table}

It is easy to see that for values of $t^*$ and values of  $N$ that correspond to them (Table 1),  the following equality holds 
 $$
 \int\limits_0^{t^*} K_N(s)ds=
\left\{ {\begin{array}{l}
 { \; \frac{1}{2},\;  \text{if}\; \, N \;\,  \text{ is odd},} \\

 { -\frac{1}{2},\;   \text{if}\;\, N  \;\, \text{is even}.} \\ 
  \end{array}}\  \right.
  $$
 Moreover, for any  $N$ 
$$
\int\limits_0^\infty K_N(s)ds=
\left\{ {\begin{array}{l}
  { 1,\;  \text{if}\; \, N \;\,  \text{is odd},} \\
 { 0,\;   \text{if}\;\, N  \;\, \text{is even},} \\
  \end{array}}\  \right.\;\; 
 $$
holds true, since  
 $$
 \left( {
- 1} \right)^{N + 1}\pi ^{2}N^{2}\int\limits_0^\infty  { e^{ - \pi ^2 N^2 { s}}}ds=
\left\{ {\begin{array}{l}
 {\;\;\, 1,\;  \text{if}\; \, N \;\,  \text{is odd},} \\
 { -1,\;   \text{if}\;\, N  \;\, \text{is even}.} \\
  \end{array}}\  \right.\;\; 
 $$
 \vspace{0.2cm}
 \begin{figure}[htbp]
   \centering
 \includegraphics[scale=0.77, trim=5cm 20.5cm 4cm 3cm]{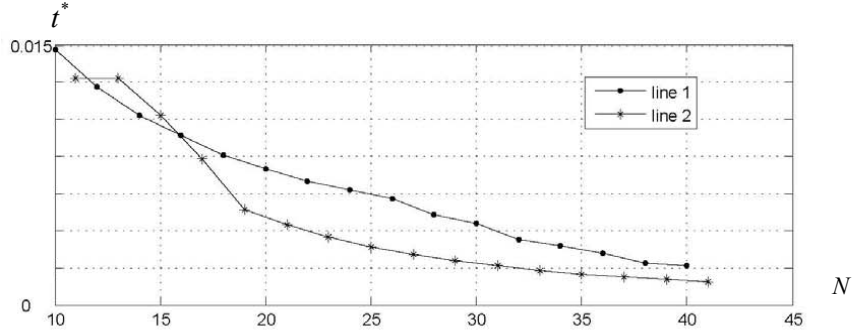}
    \caption{  "line 1"\, ---  if $N$ is even, "line 2"\, --- if $N$ is odd.} 
\end{figure}

Introduce a uniform mesh  $t_i=ih,$ $ t_{i-\frac{1}{2}}=(i-\frac{1}{2})h,$ $i=\overline{1,n},$ $ nh=T.$  
 Table 1 and Figure 1 show that the parameter $N$ from (\ref{solod6}) affects the choice of the mesh spacing  $h$.
Consider for sufficiently small  $h$ the algorithms for numerically solving the Volterra equations of the first kind   (\ref{solod5}),  (\ref{solod6}) which are based on the self-regularizing feature  of the discretization procedure \cite{Apartsyn1}, \cite{Apartsyn2}.  The midpoint rule and product integration method are used as "basic"\, ones. The product integration method is particularly effective, if  $K(t)$  is a strongly oscillating function \cite{linz1971}. Such a choice can be explained by the simplicity of the algorithm and calculation of an approximated solution with an error of order ${\cal O}(h^2)$, in the case the initial data are not perturbed.  

Solve (\ref{solod5}), (\ref{solod6}) by the midpoint rule. Approximate the integral in  (\ref{solod5}) by the sum. To calculate  $\phi(t)$, we will obtain the equation with respect to $\phi^h_{{i-\frac{1}{2}}}$ at the  $(i-\frac{1}{2})$-th node.  Denote the mesh function of interest  by  
${\phi}^h(t_{{i-\frac{1}{2}}})\equiv\check{\phi}^h_{{i-\frac{1}{2}}}$:
\begin{gather}\label{solod7}
\check{\phi}^h_{{i-\frac{1}{2}}}=\frac{y^h_i-\pi ^{2}h\sum\limits_{p = 1}^{N}\left( {
- 1} \right)^{p + 1}p^{2}\sum\limits_{j=1}^{i-1} \,{ {e^{ - \pi ^2 p^2 h {(i-j+\frac{1}{2})}}} \,\check{\phi}^h_{{j-\frac{1}{2}}}}}{\pi ^{2}h\sum\limits_{p = 1}^{N}\left( {- 1} \right)^{p + 1}p^{2} { {e^{ -\frac{1}{2} \pi ^2 p^2h}}}},
\end{gather}
 \noindent
where $i=\overline{1,n}.$

Using the product integration method, we will proceed from   (\ref{solod5}), (\ref{solod6}) to the construction  
\begin{gather}\label{solod8}
\pi ^{2}\sum\limits_{p = 1}^{N}\left( {
- 1} \right)^{p + 1}p^{2}\sum\limits_{j=1}^{i} \,\phi^h_{j-\frac{1}{2}}\, {\int\limits_{(j-1)h}^{jh} {e^{ - \pi ^2 p^2{(t- s)}}}ds} = y^h_i,\;\;i=\overline{1,n}.
 \end{gather}
By grouping the terms in (\ref{solod8}), we will write the calculation scheme for   $\phi^h(t_{{i-\frac{1}{2}}})\equiv \hat{\phi}^h_{{i-\frac{1}{2}}}$:
\begin{gather}\label{solod9}
 \hat{\phi}^h_{{i-\frac{1}{2}}}=\frac{y^h_i-\pi ^{2}\sum\limits_{p = 1}^{N}\left( {
- 1} \right)^{p + 1}p^{2}\sum\limits_{j=1}^{i-1} \,\hat{\phi}^h_{{j-\frac{1}{2}}}\,{\int\limits_{(j-1)h}^{jh} {e^{ - \pi ^2 p^2{(t- s)}}}ds }}{\pi ^{2}\sum\limits_{p = 1}^{N}\left( {- 1} \right)^{p + 1}p^{2} {\int\limits_{(i-1)h}^{ih} {e^{ - \pi ^2 p^2{(t- s)}}}ds}},\;
 \end{gather}
 where $i=\overline{1,n}.$
 
\section{Computational experiment results}

Computational experiments were made with single precision. Let us specify  $N=2,\,4$ so, that $K_{2}(0.0468)=0$, $K_{4}(0.0292)=0$.
The calculation results for the case studies from  \cite{yaparova2013} are: 

\begin{enumerate}
    \item  $\overline{\phi}_1(t)=t e^{-t}$, $t \in [0,1],$
    \item  $\overline{\phi}_2(t)= e^{-t}\sin 10\pi t$, $t \in [0,1].$
\end{enumerate}

When calculating  $\phi^h_{{i-\frac{1}{2}}}$ with respect to  (\ref{solod7}) and (\ref{solod9}),  we use a mesh analog of the function 
$$
y(t)=\pi ^{2} { \sum\limits_{p = 1}^{N} {\left( {
- 1} \right)^{p + 1} p^{2} \int\limits_{0}^{t}e^{ - \pi ^2 p^2\left( {t - s}
\right)}}\overline{\phi} \left( {s} \right)ds}$$
for fixed values of  $N=2,\,4$ and precisely specified   $\overline{\phi} \left( {s} \right)$.

Tables 2 and 3 present the values of errors $$\Vert \varepsilon^{h}_1\Vert_{C_h}=\max\limits_{1\leq i\leq n} |\overline{\phi}(t_{i-\frac{1}{2}})-\check{\phi}^h(t_{{i-\frac{1}{2}}})|$$ 
\noindent and 
$$\Vert \varepsilon^{h}_2\Vert_{C_h}=\max\limits_{1\leq i\leq n} |\overline{\phi}(t_{i-\frac{1}{2}})-\hat{\phi}^h(t_{{i-\frac{1}{2}}})|,$$
\noindent
 which are obtained using the midpoint rule and product integration method, respectively.
\begin{table}[htb!]
\caption{Errors of the mesh solution for the function  $\overline{\phi}_1$.}
\vspace{0.2cm}
 \begin{center}
\begin{tabular}{|c|c|c|c|c|}\hline
$h$&$||\varepsilon_1||_{C^h}^{N=2}$&$||\varepsilon_2||_{C^h}^{N=2}$&$||\varepsilon_1||_{C^h}^{N=4}$&$||\varepsilon_2||_{C^h}^{N=4}$ \\  \hline 
 $ 1/64$   & 0.068768&0.002936&2.957998& 0.008284 \\ 
$ 1/128$   & 0.015911&0.000734&0.024131 &  0.002235 \\
$ 1/256$   & 0.003908&0.000184&0.048312& 0.000570 \\ 
 $ 1/512$   & 0.000973&0.000046&0.011468&  0.000143 \\ 
$ 1/1024$   & 0.000243&0.000011& 0.002831&  0.000036 \\ \hline
 \end{tabular}
\end{center}
\end{table}
  \begin{table}[htb!]
\caption{Errors of the mesh solution for the function  $\overline{\phi}_2$.}
\vspace{0.2cm}
 \begin{center}
\begin{tabular}{|c|c|c|c|c|}\hline
$h$&$||\varepsilon_1||_{C^h}^{N=2}$&$||\varepsilon_2||_{C^h}^{N=2}$&$||\varepsilon_1||_{C^h}^{N=4}$&$||\varepsilon_2||_{C^h}^{N=4}$ \\  \hline 
 $ 1/64$   & 0.036243&0.028743& 1.269215&0.101544 \\ 
$ 1/128$   & 0.009016&0.007529& 0.096495&0.028924 \\
$ 1/256$   & 0.002246&0.001911& 0.023435&0.007432 \\ 
 $ 1/512$  & 0.000561&0.000481& 0.005814&0.001868 \\ 
$ 1/1024$  & 0.000140&0.000120& 0.001451&0.000468 \\ \hline
 \end{tabular}
\end{center}
\end{table}
 The Tables show that both difference methods have the convergence order ${\cal{O}}(h^2)$.

To illustrate the self-regularizing effect of the discretization procedure we will set a saw-tooth perturbance of the right-hand side of  (\ref{solod5}):
$$\tilde {y}( t_{i}) = y ( t_{i}) + (  - 1)^{i}\,\delta ,
\;
i = \overline {1,n},\; nh = T.
$$
\noindent
The Table 4 presents the values  $h_{1\, opt} (\delta)$, $h_{2\, opt} (\delta)$, which, under the fixed $\delta $, minimize the value  $$\Vert {\tilde {\varepsilon}_r }^{h\left( {\delta}  \right)} \Vert_{C_{h}}  = \mathop 
{\max}\limits_{1 \leq i \leq n( \delta)} \bigl|\mathop 
{\overline{\phi}}_r\nolimits \bigl( {t_{i - \frac{1}{2}}}  \bigr) - \tilde 
\phi_{r}^{h}(t_{i - \frac{{1}}{{2}}}) \bigr|,\; r=1,2.$$
 \noindent
    The mesh spacing was optimized by the Fibonacci method in 10 iterations. The values $\tilde \phi_{1}^h(t_{i - \frac{{1}}{{2}}})$ and $\tilde \phi^{h}_{2}(t_{i - \frac{{1}}{{2}}}) $ were calculated  by the equation    (\ref{solod7}), where $N=4$,  $T=0.0292$.
\begin{table}[ht!]
\caption{The optimal values  $h_{1} (\delta)$ and $h_{2} (\delta)$.}
\vspace{4mm}
\begin{center}
\begin{tabular}{|c|c|c|c|c|}\hline
$\delta$& $h_{1\, opt} (\delta)$ & $\Vert{\tilde {\varepsilon}}_{1}^{h_{opt}(\delta)}\Vert_{C_{h}}$ & $h_{2\, opt} (\delta)$&
$\Vert {\tilde {\varepsilon} }_{2}^{h_{opt} (\delta)} \Vert_{C_{h}}$ \\ \hline
$10^{-1}$&0.011483&0.295398&0.009186&0.332801\\ \hline
$10^{-2}$&0.009843&0.030797&0.009514&0.147060\\ \hline
$10^{-4}$&0.002297&0.001929&0.000656&0.008799\\ \hline
$10^{-5}$&0.000656&0.000629&0.000328&0.005709\\ \hline
\end{tabular}
\end{center}
\end{table}

The Table  shows that $$h_{r\,opt} (\delta) \asymp \delta 
^{\frac{{1}}{{3}}},\;\Vert {\tilde {\varepsilon}}_{r}^{h_{opt}(\delta)} \Vert_{C_{h}}  \asymp \delta 
^{\frac{{2}}{{3}}},\; r=1,2.$$
 \noindent
  Similar results were obtained when $\tilde \phi_{r}^h(t_{i - \frac{{1}}{{2}}}),\, r=1,2$ was calculated by   (\ref{solod9}).

\section{Conclusion} 

The paper considers the inverse boundary-value problem of heat conduction with a constant boundary. The problem is solved by the approach based on the direct and inverse Laplace transform. This made it possible to obtain the Volterra equation of the first kind of a special form, which characterizes an explicit relationship between the desired boundary function and the initial data on the other boundary. The algorithms have been developed to numerically solve the respective integral equations, on the basis of the midpoint rule and product integration method. The parameters determining the discretization interval were identified. The series of test calculations were made. The computational experiment shows that the numerical methods have the second order of convergence with respect to mesh spacing, and are self-regularizing under the perturbed initial data in metric  $C$.


\begin  {thebibliography}{9}


    
\bibitem {lavrentiev1980} 
Lavrentiev M.M., Romanov V.G.,  Shishatsky S.P. Ill-posed problems of mathematical physics and analysis. M.: Nauka. 1980. 287~p. (in Russian)

\bibitem {yaparova2013}  Yaparova N.M.  Numerical Methods for Solving a Boundary Value Inverse Heat Conduction
Problem// Inverse Problems in Science and Engineering, 2013,
www.tandfonline.com/doi/abs/10.1080/17415977.2013.830614. 

\bibitem {solodusha2014}  Solodusha S.V.  A numerical method for solving an inverse boundary value problem of heat conduction using the Volterra equations of the first kind// Abstract of International Triannual School-Seminar "Methods of Optimization and Their Applications", 2014, www.sei.irk.ru/conferences/mopt2014/Abstracts/Solod$\underline{~}$eng.pdf.

\bibitem{Apartsyn1}
Apartsyn A.S.,  Bakushinsky A.B. Approximated solution of the integral Volterra equations of the first kind by the method of quadrature sums // Diff.and Integ. Uravn. Irkutsk, Irkutsk State University. 1972. Issue I. P. 248--258. (in Russian)
 
\bibitem{Apartsyn2} 	
Apartsyn A.S. Discretization methods for regularization of some integral equations of the first kind // Metody chislennogo analysa i optimizatsii. Novosibirsk: Nauka, Sib. Branch. 1987. P. 263--297. (in Russian)

\bibitem {linz1971}
Linz P. Product integration method for Volterra integral equations of the first kind // BIT, 1971.  Vol. 11. P. 413--421.



\end{thebibliography}
\end{document}